\documentclass[12pt]{amsart}
\usepackage{amssymb,eepic,verbatim}

\usepackage{pstricks,pst-node}
\psset{unit=.5cm}

\newcommand{\C}{\mathbb{C}}
\renewcommand{\P}{\mathbb{P}}
\newcommand{\Z}{\mathbb{Z}}
\newcommand{\on}{\operatorname}
\newcommand{\Gr}{ \on{Gr} }
\newcommand{\Fl}{ \on{Fl} }
\renewcommand{\.}{ {\scriptscriptstyle{\bullet}}  }
\renewcommand{\min}{{\text{min}}}
\renewcommand{\max}{{\text{max}}}

\begin{document}

\title{On the quantum cohomology of homogeneous varieties}

\date{February 27,2003}

\author{W. Fulton}\thanks{Partially supported by NSF grant
DMS9970435}
\address{Department of Mathematics,
University of Michigan,
Ann Arbor, Michigan  48109-1109, U.S.A.}
\email{wfulton@umich.edu}

\maketitle

About a decade ago physicists set off something of a ``big bang''
in
the universe of algebraic geometry. A new approach to enumerative
geometry, involving mirror symmetry, solved some old questions and raised
many new ones. To  mention just a few of the areas influenced by these
developments, we have new approaches to the study of Calabi-Yau manifolds
and orbifolds, the notions of Gromov-Witten and related invariants, and
new results about intersection theory on moduli spaces.

From this revolution we now have general and powerful theorems, which
produce meaningful Gromov-Witten invariants on every smooth projective
complex variety (by means of virtual fundamental cycles), with methods for
calculating (localization).  On the other hand, a good general
understanding of all invariants on a given space has come only recently,
in the work first of Kontsevich in the case of a point, and the work of
Okounkov and Pandharipande in the case of the projective line.

One of the first situations explained to us by E.~Witten \cite{Witten}
was for the Grassmann variety $\Gr(k,n)$ of $k$-dimensional linear
subspaces of an $n$-dimensional complex vector space, where a deformation
of the classical cohomology, called the (small) quantum cohomology, was
described.  For this variety, and for general homogeneous varieties $G/P$,
the Gromov-Witten invariants are given by na\"{\i}ve counting --- no 
virtual
cycles are needed.  Some of this story has been extended to other
homogeneous varieties.  It is somewhat surprising that, after all the
progress in this general area, there are still so many open questions
about the quantum cohomology of $G/P$ in general, and even for $\Gr(k,n)$.
One of the appealing features of this work is the fact that the classical
study of these varieties leads to interesting combinatorics.  There is
already considerable evidence that the quantum versions involve equally
interesting combinatorial ideas.

The aim of this lecture is to sketch what is known and what remains open
about the quantum cohomology of homogeneous varieties.  We will be
concerned only with the ``small'' quantum cohomology, which is formed from
the $3$-point Gromov-Witten invariants.  Unlike much of the algebraic
geometry that has come from physics, this part of the story is quite
accessible to a general audience.  Even here, we do not attempt a survey,
and apologize to the many whose work is not cited.

We will be concerned with a variety $X$ which is a homogeneous variety: $X
= G/P$, where $G$ is a simple complex algebraic Lie group, and $P$ is a
parabolic subgroup.  These are classified by data from Dynkin diagrams,
cf. \cite{F-H}.  When $G = \operatorname{SL}_n(\C)$, the varieties $G/P$
are the varieties of partial flags of subspaces of $\C^n$, generalizing
the Grassmannian (the one-step flags).  There are corresponding manifolds
of flags in orthogonal or symplectic vector spaces for the other groups of
classical type, together with a few exceptional cases.  Those unfamiliar
with the general roots and weights story can concentrate on the
Grassmannians $\Gr(k,n)$.

The classical cohomology $H^*(X) = H^*(X,\Z)$ has a basis (over $\Z$) of
classes $\sigma_\lambda$ of Schubert varieties, where $\lambda$ varies
over a combinatorial set.  In general, this set is the quotient $W/W_P$ of
the Weyl group $W$ of $G$ by the subgroup $W_P$ corresponding to $P$.  For
$X = \Gr(k,n)$, this combinatorial set can be taken to be the set of
partitions whose Young diagram fits in a $k$ by $n-k$ rectangle, i.e.,
$\lambda = (\lambda_1, \ldots, \lambda_k)$, with $n-k \geq \lambda_1 \geq
\lambda_2 \geq \ldots \lambda_k \geq 0$.  For example, take $k = 4$, $n 
= 9$, and $\lambda = (5,4,4,3)$, with Young diagram:

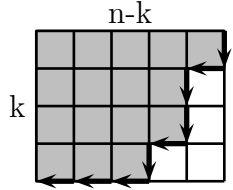
\pspicture(-9,-1)(6,5)
\pspolygon[fillstyle=solid, fillcolor=lightgray]
(0,0)(3,0)(3,1)(4,1)(4,3)(5,3)(5,4)(0,4)
\psline[linewidth=1pt](0,0)(5,0)
\psline[linewidth=1pt](0,1)(5,1)
\psline[linewidth=1pt](0,2)(5,2)
\psline[linewidth=1pt](0,3)(5,3)
\psline[linewidth=1pt](0,4)(5,4)
\psline[linewidth=1pt](0,0)(0,4)
\psline[linewidth=1pt](1,0)(1,4)
\psline[linewidth=1pt](2,0)(2,4)
\psline[linewidth=1pt](3,0)(3,4)
\psline[linewidth=1pt](4,0)(4,4)
\psline[linewidth=1pt](5,0)(5,4)
\psline[linewidth=2pt]{->}(5,4)(5,3)
\psline[linewidth=2pt]{->}(5,3)(4,3)
\psline[linewidth=2pt]{->}(4,3)(4,2)
\psline[linewidth=2pt]{->}(4,2)(4,1)
\psline[linewidth=2pt]{->}(4,1)(3,1)
\psline[linewidth=2pt]{->}(3,1)(3,0)
\psline[linewidth=2pt]{->}(3,0)(2,0)
\psline[linewidth=2pt]{->}(2,0)(1,0)
\psline[linewidth=2pt]{->}(1,0)(0,0)
\rput(-.5,2){k}
\rput(2.5,4.5){n-k}
\endpspicture

\noindent As indicated, we will often regard the partition as a path from
the upper right corner of this rectangle to the lower left corner, by a
sequence of $n$ steps, each either down or to the left.

If one fixes a complete flag $F_{\.}$ of subspaces of $\C^n$, with
$\dim{F_i} = i$ for $0 \leq i \leq n$, then $\sigma_{\lambda}$ is the
cohomology class of the Schubert variety
\[
\Omega_{\lambda} = \Omega_{\lambda}(F_{\,}) \, = \, \{L \subset \C^n \mid
 \dim(L \cap F_i) \geq r_i \text{ for } 1 \leq i \leq n\} ,
\]
where $r_i$ is the row one has passed after $i$ steps along the path
defining $\lambda$.
The (complex) codimension of $\Omega_{\lambda}$ is the number $|\lambda| 
=
\sum \lambda_i$ of boxes in the Young diagram (above the path), so
$\sigma_{\lambda}$ is in $H^{2|\lambda|}(X)$.

In the cohomology ring $H^*(X)$, one has
\[
\sigma_\lambda \cdot \sigma_\mu = \sum c_{\lambda \, \mu}^{\,\, \nu} \,
\sigma_\nu,
\]
the sum over $\nu$ with $|\nu| = |\lambda| + |\mu|$.  The coefficients are
nonnegative integers known as Littlewood-Richardson coefficients, which have
several interesting combinatorial descriptions, and play a leading role in
algebraic combinatorics.  There are dual classes
$\sigma_{\lambda^{\vee}}$, where $\lambda^{\vee} = (n-k-\lambda_k, \ldots
,
n-k-\lambda_1)$ is obtained by rotating the diagram of the complement
of $\lambda$ by $180$ degrees; the intersection number
$\int \sigma_\lambda \cdot \sigma_\mu$ is $1$ if $\mu = \lambda^{\vee}$
and $0$ otherwise.  With this notation,
\[
c_{\lambda \, \mu}^{\,\, \nu} =
\int \sigma_\lambda \cdot \sigma_\mu \cdot \sigma_{\nu^{\vee}}
\]
is the number of points in the intersection of three Schubert varieties
$\Omega_\lambda$, $\Omega_{\mu}'$, and $\Omega_{\nu^{\vee}}''$, using
three general flags.

The quantum cohomology $QH^*(X)$ has the same basis $\sigma_\lambda$, but
as an algebra over a polynomial ring $\Z[q]$ (in general, the polynomial
ring $\Z[q_\beta]$ with a variable for each $\beta$ 
in $W/W_P$ 
with $|\beta| = 1$).
Here the (complex) codimension of $q$ is $n$ (or
$\int_{\sigma_{\beta^{\vee}}} c_1(T_X)$ in general).  In this ring, with
its product denoted by a $*$,
\[
\sigma_\lambda * \sigma_\mu = \sum c_{\lambda \, \mu}^{\,\, \nu}(d) \, q^
d
\,  \sigma_\nu,
\]
the sum over 
$d$ and 
$\nu$ with $|\nu| = |\lambda| + |\mu| - dn$.  The coefficient
$ c_{\lambda \, \mu}^{\,\, \nu}(d)$ is the number of maps of degree $d$
from the projective line to $X$, whose image meets three general Schubert
varieties $\Omega_\lambda$, $\Omega_{\mu}'$, and $\Omega_{\nu^{\vee}}''$.
(In general, the degree is defined by the equation $f_*[\P^1] = \sum
d_\beta \sigma_{\beta^{\vee}}$.)  The surprising fact is that this product
defines an {\em associative} and commutative $\Z[q]$-algebra.

These coefficients (the $3$-point Gromov-Witten numbers) are known in
principle.  But many questions about them remain open --- very much so for
a general $G/P$, and considerably so for $\Gr(k,n)$.  To have a good
understanding of these numbers, one wants: (1) a presentation of the
quantum cohomology ring:
\[
QH^*(X) \, = \, \Z[q][x_1, \ldots , x_N]/(\text{relations}).
\]
Presentations for quantum cohomology have been given for a general $G/B$
by B.~Kim \cite{Kim}.  In type A explicit presentations have been given
for the partial flag varieties (see \cite{AS} and \cite{CF2}); the case of
Lagrangian and maximal orthogonal Grassmannians can be found in \cite{KT}.
The presentation for general $G/P$ has been announced in MIT lectures by
D.~Peterson (unpublished).

One also wants: (2) a ``quantum Giambelli'' formula, which expresses each
$\sigma_\lambda$ as a polynomial in the generators $x_i$ and $q$.
Such formulas were worked out first for the Grassmannian by A.~Bertram
\cite{Ber} (where the answer is the same as in the classical case ---
there is no quantum correction), and for complete and partial flag
varieties by I.~Ciocan-Fontanine \cite{CF1}, \cite{CF2}, and S.~Fomin,
S.~Gelfand, and A.~Postnikov \cite{FGP}; see L.~Chen \cite{Chen} for a
concise treatment. A.-L.~Mare \cite{Mare} has recently given an algorithm
for arbitrary $G/B$'s.

In addition, one would like: (3) a combinatorial formula for the
coefficients
$c_{\lambda \, \mu}^{\,\, \nu}(d)$, or at least to know which are nonzero.
Although (3) follows in principle from (1) and (2), it is far from obvious
how to carry this out explicitly, even for the classical cohomology.  In
fact, for $H^*(\Gr(k,n))$, there is a recent criterion to tell for which
partitions $c_{\lambda \, \mu}^{\,\, \nu}$ is not zero, in terms of a
collection of linear inequalities (see \cite{BAMS} for an exposition of
this story).  In fact, calculations done with A.~Buch indicate that there
are similar inequalities to describe the positivity of the quantum numbers
$c_{\lambda \, \mu}^{\,\, \nu}(d)$.  For the classical numbers, the
inequalities are determined by the answers to the same questions for
smaller Grassmannians $\Gr(r,k)$, $1 \leq r \leq k$.  It is natural to
hope that the same is true for the quantum cohomology.\footnote{(Added
later) P.~Belkale has done this, in ``The quantum Horn conjecture," 
math.AG/0303013.}

The proof of basic properties of quantum cohomology uses the space of
stable maps from rational curves to the given variety \cite{KM}, 
cf.~\cite{FP}.
Strangely, proofs of quantum Giambelli and other quantum formulas have
used Grothendieck's quot scheme compactifications rather than this 
space of stable maps.  Recently, however,
Buch has shown how the basic facts about quantum cohomology of $\Gr(k,n)$
can be proved be entirely elementary methods \cite{Buch1}, \cite{Buch2},
\cite{Buch3}, 
without any compactifications at all.  For a morphism $f \colon \P^1 \to
\Gr(k,n)$ of degree $d$, consider the intersection $K$ of all the linear
spaces $f(t)$ as $t$ varies in $\P^1$, and the span $S$ of all these
linear spaces.  Buch proves that
\[
\dim K \geq k-d  \quad \text{ and } \quad  \dim S \leq k+d.
\]
With A.~Kresch and H.~Tamvakis \cite{BKT}, he uses this to show that
$c_{\lambda \, \mu}^{\,\, \nu}(d)$ is equal to a classical intersection
number of Schubert classes in the flag variety
$\Fl(k-d,k+d;\C^n)$.  There is a conjectured combinatorial formula for
intersection numbers on two-step flag varieties, from
A.~Knutson, so one now has at least a conjectured combinatorial formula
for the quantum Littlewood-Richardson coefficients.

During the 1996--97 year at the Mittag-Leffler Institute, many tables of
products $\sigma_\lambda * \sigma_\mu$ were computed.  It was surprising
to see that the product was never $0$, since, for degree reasons, there
are only a finite number of possible nonzero coefficients.  This was
proved, for the Grassmannian, by S.~Agnihotri and C.~Woodward \cite{AW}.
\footnote{Agnihotri gave a very short proof, by observing that, after
multiplying by $\sigma_{\lambda^{\vee}}$ and $\sigma_{\mu^{\vee}}$ --- and
noting that by the nonnegativity of all coefficients, there can be no
cancelation --- it suffices to do this when $\lambda$ and $\mu$ are both
the maximal partitions $(n-k, \ldots, n-k)$, where it is an easy
calculation.}

Note that the classical product $\sigma_\lambda \cdot \sigma_\mu$ is not
zero exactly when the diagram of $\lambda$ is contained in the diagram of
$\mu^{\vee}$.  On the basis of some calculations\footnote{Buch has an
efficient computer program for calculating classical and quantum
Littlewood-Richardson coefficients, available at
http://home.imf.au.dk/abuch/lrcalc/.}, it was natural to conjecture that
the smallest power $q^d$ that appears in the quantum product
$\sigma_\lambda * \sigma_\mu$ would be the maximum $d$ such that a $d$ by
$d$ square fits inside the diagram of $\lambda$ but outside that of
$\mu^{\vee}$.
For example, take $k = 4$, $n = 9$, $\lambda = (5,4,4,3)$ and $\mu 
=
(5,4,4,1)$ (so $\mu^{\vee} = (4,1,1,0)$:

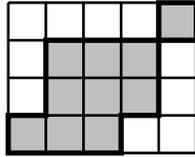
\pspicture(-9,-1)(6,4)
\pspolygon[fillstyle=solid, fillcolor=lightgray]
(0,0)(3,0)(3,1)(4,1)(4,3)(1,3)(1,1)(0,1)
\pspolygon[fillstyle=solid, fillcolor=lightgray]
(4,3)(5,3)(5,4)(4,4)
\psline[linewidth=1pt](0,0)(5,0)
\psline[linewidth=1pt](0,1)(5,1)
\psline[linewidth=1pt](0,2)(5,2)
\psline[linewidth=1pt](0,3)(5,3)
\psline[linewidth=1pt](0,4)(5,4)
\psline[linewidth=1pt](0,0)(0,4)
\psline[linewidth=1pt](1,0)(1,4)
\psline[linewidth=1pt](2,0)(2,4)
\psline[linewidth=1pt](3,0)(3,4)
\psline[linewidth=1pt](4,0)(4,4)
\psline[linewidth=1pt](5,0)(5,4)
\pspolygon[linewidth=2pt]
(0,0)(3,0)(3,1)(4,1)(4,3)(5,3)(5,4)(4,4)(4,3)
(1,3)(1,1)(0,1)
\endpspicture

This indicates that $d = 2$, as there are three $2$ by $2$ squares that
can be inserted inside $\lambda$ and outside $\mu^{\vee}$.  In fact,
\[
\sigma_{(5,4,4,3)} * \sigma_{(5,4,4,1)} =
q^2(\sigma_{(5,3,2,2)} + \sigma_{(5,3,3,1)} + \sigma_{(5,4,2,1)}) +
q^3(\sigma_{(3,0,0,0)} + \sigma_{(2,1,0,0)} + \sigma_{(1,1,1,0)}).
\]
Woodward and I \cite{FW} have proved this conjecture.  Moreover, the
product is nonzero for every $G/P$, and there is a combinatorial formula
for the minimum $d$ such that $q^d$ occurs in the product.

For the Grassmannian, other proofs have been given by Buch \cite{Buch1}
and Belkale \cite{Bel}, who gives the following refinement.  The
coefficient of the minimal $q^d$ is equal to a classical product
$\sigma_{\lambda'} \cdot \sigma_{\mu'}$, where $\lambda'$ and $\mu'$ are
obtained as follows.  Pick any maximal square inside $\lambda$ and outside
$\mu^{\vee}$, and let $a$ be the number of steps one travels down the path
of $\lambda$ to reach the southeast corner of the square, and let $b$ be
the number of steps one travels up on the path of $\mu^{\vee}$ to reach
the northwest corner of the square.

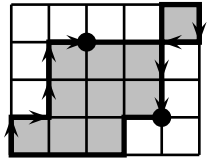
\pspicture(-9,-1)(6,5)
\pspolygon[fillstyle=solid, fillcolor=lightgray]
(0,0)(3,0)(3,1)(4,1)(4,3)(1,3)(1,1)(0,1)
\pspolygon[fillstyle=solid, fillcolor=lightgray]
(4,3)(5,3)(5,4)(4,4)
\psline[linewidth=1pt](0,0)(5,0)
\psline[linewidth=1pt](0,1)(5,1)
\psline[linewidth=1pt](0,2)(5,2)
\psline[linewidth=1pt](0,3)(5,3)
\psline[linewidth=1pt](0,4)(5,4)
\psline[linewidth=1pt](0,0)(0,4)
\psline[linewidth=1pt](1,0)(1,4)
\psline[linewidth=1pt](2,0)(2,4)
\psline[linewidth=1pt](3,0)(3,4)
\psline[linewidth=1pt](4,0)(4,4)
\psline[linewidth=1pt](5,0)(5,4)
\pspolygon[linewidth=2pt]
(0,0)(3,0)(3,1)(4,1)(4,3)(5,3)(5,4)(4,4)(4,3)
(1,3)(1,1)(0,1)
\psline[linewidth=2pt]{->}(5,4)(5,3)
\psline[linewidth=2pt]{->}(5,3)(4,3)
\psline[linewidth=2pt]{->}(4,3)(4,2)
\psline[linewidth=2pt]{->}(4,2)(4,1)
\psline[linewidth=2pt]{->}(0,0)(0,1)
\psline[linewidth=2pt]{->}(0,1)(1,1)
\psline[linewidth=2pt]{->}(1,1)(1,2)
\psline[linewidth=2pt]{->}(1,2)(1,3)
\psline[linewidth=2pt]{->}(1,3)(2,3)
\psdots*[dotscale=2](4,1)(2,3)
\endpspicture

In this example, these corners and paths are marked, so $a = 4$ and $b =
5$.
Add a  $k$ by $a$ rectangle to the left of $\lambda$, and remove $n$-rims
until the result lies inside the original $k$ by $n-k$ rectangle; this
result is $\lambda'$.  Similarly, to find $\mu'$, one removes rims from
the addition of a $k$ by $b$ rectangle to $\mu$.  In the example:

\pspicture(-1,-1)(23,5)
\pspolygon[fillstyle=solid, fillcolor=lightgray]
(0,2)(1,2)(1,3)(4,3)(4,4)(0,4)
\psline[linewidth=1pt](0,0)(7,0)
\psline[linewidth=1pt](0,1)(8,1)
\psline[linewidth=1pt](0,2)(8,2)
\psline[linewidth=1pt](0,3)(9,3)
\psline[linewidth=1pt](0,4)(9,4)
\psline[linewidth=1pt](0,0)(0,4)
\psline[linewidth=1pt](1,0)(1,4)
\psline[linewidth=1pt](2,0)(2,4)
\psline[linewidth=1pt](3,0)(3,4)
\psline[linewidth=1pt](4,0)(4,4)
\psline[linewidth=1pt](5,0)(5,4)
\psline[linewidth=1pt](6,0)(6,4)
\psline[linewidth=1pt](7,0)(7,4)
\psline[linewidth=1pt](8,1)(8,4)
\psline[linewidth=1pt](9,3)(9,4)
\psline[linewidth=2pt, linecolor=gray]
(.5,.5)(.5,1.5)(1.5,1.5)(1.5,2.5)(4.5,2.5)
(4.5,3.5)(5.5,3.5)
\psline[linewidth=2pt, linecolor=gray]
(1.5,.5)(2.5,.5)(2.5,1.5)(5.5,1.5)(5.5,2.5)
(6.5,2.5)(6.5,3.5)
\psline[linewidth=2pt, linecolor=gray]
(3.5,.5)(6.5,.5)(6.5,1.5)(7.5,1.5)(7.5,3.5)
(8.5,3.5)
\pspolygon[fillstyle=solid, fillcolor=lightgray]
(12,0)(13,0)(13,2)(14,2)(14,3)(15,3)(15,4)(12,4)
\psline[linewidth=1pt](12,0)(18,0)
\psline[linewidth=1pt](12,1)(21,1)
\psline[linewidth=1pt](12,2)(21,2)
\psline[linewidth=1pt](12,3)(22,3)
\psline[linewidth=1pt](12,4)(22,4)
\psline[linewidth=1pt](12,0)(12,4)
\psline[linewidth=1pt](13,0)(13,4)
\psline[linewidth=1pt](14,0)(14,4)
\psline[linewidth=1pt](15,0)(15,4)
\psline[linewidth=1pt](16,0)(16,4)
\psline[linewidth=1pt](17,0)(17,4)
\psline[linewidth=1pt](18,0)(18,4)
\psline[linewidth=1pt](19,1)(19,4)
\psline[linewidth=1pt](20,1)(20,4)
\psline[linewidth=1pt](21,1)(21,4)
\psline[linewidth=1pt](22,3)(22,4)
\psline[linewidth=2pt, linecolor=gray] (13.5,.5)(13.5,1.5)
(14.5,1.5)(14.5,2.5)(15.5,2.5)(15.5,3.5)
(18.5,3.5)
\psline[linewidth=2pt, linecolor=gray] (14.5,.5)(15.5,.5)
(15.5,1.5)(16.5,1.5)(16.5,2.5)(19.5,2.5)
(19.5,3.5)
\psline[linewidth=2pt, linecolor=gray] (16.5,.5)(17.5,.5)
(17.5,1.5)(20.5,1.5)(20.5,3.5)(21.5,3.5)
\endpspicture

One sees that $\lambda' = (4,1,0,0)$ and $\mu' = (3,2,1,1)$.  And, as
predicted, the classical product of $\sigma_{(4,1,0,0)}$ and
$\sigma_{(3,2,1,1)}$ is $\sigma_{(5,3,2,2)} + \sigma_{(5,3,3,1)} +
\sigma_{(5,4,2,1)}$.
In general, $\sigma_{(a, \ldots, a)} * \sigma_\lambda = q^m
\sigma_{\lambda'}$, where $m$ is the number of $n$-rims removed; this is a
special case of the algorithm of \cite{BCF} for arbitrary multiplication
of quantum Schubert classes.

The general formula for the minimal $d$ such that $q^d$ occurs involves
the combinatorics of the curves in $G/P$ that are invariant by the action
of a maximal torus.  It is measured by the length of a chain of such
curves needed to join a point of a Schubert variety for $\lambda$ to an
opposite Schubert variety for $\mu^{\vee}$.

On the Grassmannian, A.~Yong \cite{Yong} found an upper bound for which
powers of $q$ can occur, and Postnikov \cite{Post2} has found exactly
which $d$ have $q^d$ occur in a product $\sigma_\lambda * \sigma_\mu$.
They are all $d$ between the above $d_\min$ and a $d_\max$ determined as
follows.  Slide the path for the diagram of $\mu^{\vee}$ southeast
$d_\min$ steps, so there is nothing outside it and inside the path for
$\lambda$.  The ends of this path now extend outside the given $k$ by
$n-k$ rectangle.  Think of this rectangle as a torus, with opposite sides
identified, and regard the extended path as a loop on the torus.

\pspicture(-8,-3)(8,5)
\psline[linewidth=1pt](0,0)(5,0)
\psline[linewidth=1pt](0,1)(5,1)
\psline[linewidth=1pt](0,2)(5,2)
\psline[linewidth=1pt](0,3)(5,3)
\psline[linewidth=1pt](0,4)(5,4)
\psline[linewidth=1pt](0,0)(0,4)
\psline[linewidth=1pt](1,0)(1,4)
\psline[linewidth=1pt](2,0)(2,4)
\psline[linewidth=1pt](3,0)(3,4)
\psline[linewidth=1pt](4,0)(4,4)
\psline[linewidth=1pt](5,0)(5,4)
\psline[linewidth=3pt, linecolor=gray]
(2,-2)(2,-1)(3,-1)(3,0)
\psline[linewidth=3pt, linecolor=gray]
(5,1)(6,1)(6,2)(7,2)
\psline[linewidth=3pt]
(3,0)(3,1)(5,1)
\psline[linewidth=3pt]
(0,1)(1,1)(1,2)(2,2)(2,3)(3,3)(3,4)
\psline[linewidth=2pt]
(4,1)(4,3)(5,3)(5,4)
\endpspicture

\noindent Keep moving this loop to the southeast until it just touches the
path for $\lambda$ again.  The total number of steps is $d_\max$.  In this
example, $d_\max = 3$.  Postnikov has a kind of duality between the
minimal and maximal terms, 
including an explicit formula for the maximal terms, 
although those in between are not yet well
understood.  He also shows that, for all $G/B$, only one minimal degree
occurs.

The proof in \cite{FW} uses Kleiman's transversality theorem on the
Kontsevich-Manin moduli spaces, and the fixed points of the torus action.
The work in \cite{BCF} is combinatorial, as is that in \cite{Post2}.

We conclude with a few remarks and questions.

Postnikov \cite{Post2} gives a formula for the quantum
Littlewood-Richardson numbers in terms of ``tableaux on a torus'', but not
in a way to show their positivity.

The quantum cohomology of $\Gr(k,n)$ is related to the fusion ring
(Verlinde algebra) constructed from representations of the unitary group
$U_k$ at level $n-k$, with the Schubert class $\sigma_\lambda$
corresponding to the representation $V_\lambda$ with highest weight
$\lambda$.  The fact that $V_\lambda \otimes V_\mu$ is not zero gives
another proof that the quantum product $\sigma_\lambda * \sigma_\mu$ is
not zero.  There are many physics papers studying these tensor products.
G.~Tudose \cite{Tud} has a formula for the product when 
one partition has two parts. 

Classically, knowing $H^*(G/B)$ is essentially equivalent to knowing
$H^*(G/P)$ for all $P$, or all maximal parabolic $P$.  In the quantum
world, this is far from the case, because of the lack of functoriality of
quantum cohomology.

The fact that the coefficients $c_{\lambda \, \mu}^{\,\, \nu}(d)$ are
nonnegative puts very strong restrictions on the quantum cohomology rings.
In practise, at least in small examples, one can often use that fact,
together with very few simple calculations, to determine the whole quantum
cohomology ring.  It would be desirable to prove theorems characterizing
quantum cohomology along these lines.

Peterson has given a remarkable formula for each $c_{\lambda \, \mu}^{\,\,
\nu}(d)$ on a $G/P$ as a number $c_{\lambda' \, \mu'}^{\,\, \nu'}(d')$ on
the corresponding $G/B$.  This is explained in \cite{Wood}.

Various symmetries have been found in the quantum cohomology of
Grassmannians.  For these, see \cite{AW}, \cite{Post1}, \cite{Bel},
\cite{Heng}.


\begin{thebibliography}{99}

\bibitem{AW}
S.~Agnihotri and C.~Woodward,
\newblock Eigenvalues of products of unitary matrices and quantum
{S}chubert calculus,
\newblock {\em Math. Res. Lett.}, 5(6):817--836, 1998,
\newblock alg-geom/9712013.

\bibitem{AS}
A.~Astashkevich and V. Sadov,
\newblock Quantum cohomology of partial flag manifolds ${F}\sb {n\sb
  1\cdots n\sb k}$,
\newblock {\em Comm. Math. Phys.}, 170:503--528,1995,
\newblock hep-th/9401103.


\bibitem{Bel}
P.~Belkale,
\newblock Transformation formulas in quantum cohomology,
\newblock math.AG/0210050.

\bibitem{Ber}
A.~Bertram,
\newblock Quantum schubert calculus,
\newblock {\em Adv. Math.}, 128:289--305, 1997,
\newblock alg-geom/9410024.

\bibitem{BCF}
A.~Bertram, I.~Ciocan-Fontanine, and W.~Fulton,
\newblock Quantum multiplication of {S}chur polynomials,
\newblock {\em J. Algebra}, 219(2):728--746, 1999,
\newblock alg-geom/9705024.

\bibitem{Buch1} A.~Buch,
\newblock Quantum cohomology of {G}rassmannians,
\newblock math.AG/0106268.

\bibitem{Buch2} A.~Buch,
\newblock Direct proof of the quantum Monk's formula,
\newblock {\em Proc. Amer. Math. Soc.}, 131(7):2037--2042, 2003,
\newblock math.AG/0107202.

\bibitem{Buch3} A.~Buch,
\newblock Quantum cohomology of partial flag manifolds,
\newblock math.AG/0303245.

\bibitem{BKT}
A.~Buch, A.~Kresch, and H.~Tamvakis,
\newblock Gromov-Witten invariants on {G}rassmannians,
\newblock math.AG/0306388.

\bibitem{Chen}
L.~Chen,
\newblock Quantum cohomology of flag manifolds,
\newblock {\em Adv. Math.} 174:1--34, 2003,
\newblock math.AG/0010080.

\bibitem{CF1}
I.~Ciocan-Fontanine,
\newblock The quantum cohomology ring of flag varieties,
\newblock {\em Trans. Amer. Math. Soc.}, 351(7):2695--2729, 1999,
\newblock alg-geom/9505002.

\bibitem{CF2}
I.~Ciocan-Fontanine,
\newblock On quantum cohomology rings of partial flag varieties,
\newblock {\em Duke Math. J.}, 98(3):485--524, 1999,
\newblock math.AG/9710213.

\bibitem{FGP}
S.~Fomin, S.~Gelfand, and A.~Postnikov,
\newblock Quantum {S}chubert polynomials,
\newblock {\em J. Amer. Math. Soc.}, 10:565--596, 1997.

\bibitem{BAMS}
W.~Fulton,
\newblock Eigenvalues, invariant factors, highest weights, and {S}chubert
calculus,
\newblock {\em Bull. Amer. Math. Soc.}, 37:209--249, 2000,
\newblock math.AG/0301307.

\bibitem {F-H}
W.~Fulton and J.~Harris,
\newblock {\em Representation {T}heory},
\newblock Springer-Verlag, 1991, 1999.

\bibitem{FP}
W.~Fulton and R.~Pandharipande,
\newblock Notes on stable maps and quantum cohomology,
\newblock In {\em Algebraic geometry---Santa Cruz 1995}, pages 45--96.
Amer. Math. Soc., 1997,
\newblock alg-geom/9608011.

\bibitem{FW}
W.~Fulton and C.~Woodward,
\newblock On the quantum product of {S}chubert classes,
\newblock {\em J. Alg.Geom.}, to appear,
\newblock math.AG/0112183.

\bibitem{Heng}
H.~Hengelbrock,
\newblock An involution on the quantum cohomology ring of the
Grassmannian,
\newblock math.AG/0205260.

\bibitem{Kim}
Bumsig Kim,
\newblock Quantum cohomology of flag manifolds ${G}/{B}$ and quantum
{T}oda lattices,
\newblock {\em Ann. of Math. (2)}, 149(1):129--148, 1999,
\newblock alg-geom/9607001.

\bibitem{KM}
M.~Kontsevich and Yu.~Manin,
\newblock Gromov-{W}itten classes, quantum cohomology, and
enumerative geometry,
\newblock{\em Comm. Math. Phys.}, 164(3):525--562, 1994,
\newblock hep-th/9402147.

\bibitem{KT}
A.~Kresch and H.~Tamvakis,
\newblock Quantum cohomology of the {L}agrangian {G}rassmannian,
\newblock math.AG/0306337.

\bibitem {Mare}
A.-L.~Mare,
\newblock  Polynomial representatives of Schubert classes in
${QH}^*({G}/{B})$,
\newblock {\em Math. Res. Lett.} 9(5--6):757--760, 2002,
\newblock math.CO/0205309.


\bibitem{Post1}
A.~Postnikov,
\newblock Symmetries of {G}romov-{W}itten invariants,
\newblock in {\em Advances in algebraic geometry motivated by physics
(Lowell, MA, 2000)}, pages 251--258,
{\em Contemp. Math.} 276, Amer. Math. Soc., 2001,
\newblock math.CO/0009174.

\bibitem{Post2}
Alexander Postnikov,
\newblock Affine approach to quantum {S}chubert calculus,
\newblock math.CO/0205165.

\bibitem {Tud}
G.~Tudose,
\newblock A special case of sl(n)-fusion coefficients,
\newblock math.CO/0008034.

\bibitem {Witten}
E.~Witten,
\newblock The {V}erlinde algebra and the cohomology of the {G}rassmannian,
\newblock in {\em Geometry, topology, \& physics}, Conf. Proc. Lecture
Notes Geom. Topology, IV, pages 357--422, Internat. Press, 1995,
\newblock hep-th/9312104.

\bibitem{Wood}
C.~Woodward,
\newblock On {D}. {P}eterson's comparison formula for {G}romov-{W}itten
  invariants of ${G}/{P}$,
\newblock math.AG/0206073.

\bibitem{Yong}
A.~Yong,
\newblock Degree bounds in quantum Schubert calculus,
\newblock {\em Proc. Amer. Math. Soc.} 141(9):2649--2655, 2003,
\newblock math.CO/0112133.

\end{thebibliography}
\end{document}